\documentclass{amsart}

\newtheorem{thm}{Theorem}[section]
\newtheorem{lem}[thm]{Lemma}
\newtheorem{cor}[thm]{Corollary}

\theoremstyle{definition}
\newtheorem{defn}[thm]{Definition}
\newtheorem{ex}[thm]{Example}

\newtheorem{rem}[thm]{Remark}

\newcommand{\point}{\mbox {\rm point}}

\newcommand{\Z}{\mathbf Z}
\newcommand{\Q}{\mathbf Q}
\newcommand{\Ql}{{\Q}_{\ell}}

\newcommand{\PP}{\mathbf P}
\newcommand{\Le}{\mathfrak l}
\newcommand{\Det}{\mbox{\rm det}}

\newcommand{\Num}{\mbox{\rm Num}}
\newcommand{\Al}{\mbox{\rm Al}}
\newcommand{\Aut}{\mbox{\rm Aut}}
\newcommand{\tr}{\mbox{\rm tr}}

\newcommand{\Graph}{\mbox{\rm Graph}}
\newcommand{\cl}{\mbox{\rm cl}}

\newcommand{\id}{\mbox{\rm id}}

\newcommand{\Hom}{\mbox{\rm Hom}}
\newcommand{\End}{\mbox{\rm  End}}

\newcommand{\Frob}{\mbox{\rm Frob}}

\begin{document}
\title[Powers of ordinary cubic fourfolds]{The Tate Conjecture for Powers of Ordinary Cubic Fourfolds Over Finite Fields}
\author {Yuri G. Zarhin}
\address{The Pennsylvania State University, Department of Mathematics\\
\noindent 218 McAllister Building, University Park, PA 16802, USA}
\email{zarhin@@math.psu.edu}
\date{}
\thanks{Partially supported by the NSF}
\begin{abstract}
Recently N. Levin proved the Tate conjecture for ordinary cubic
fourfolds over finite fields. In this paper we prove the Tate
conjecture for self-products of ordinary cubic fourfolds. Our
proof is based on properties of so called polynomials of K3 type
introduced by the author about a dozen years ago.
\end{abstract}
\maketitle

\section{Introduction.}
Let $X$ be a smooth projective variety over a finite field $k$ of
characteristic $p$. We write $X_a$ for $X \times k(a)$ where
$k(a)$ is an algebraic closure of $k$. For each non-negative
integer $m$ and each rational prime $l$ different from $p$ let us
consider the $2m$th twisted $\ell$-adic cohomology group
$H^{2m}(X_a,\Ql)(m)$ of $Y_a$. The Galois group $G(k)$ of $k$ acts
on  $H^{2m}(X_a,\Ql)(m)$ in a natural way. In \cite{T1} Tate
conjectured that the subspace fixed under the Galois action is
spanned by cohomology classes of codimension $m$ algebraic cycles
on $X$. The famous conjecture of Serre and Grothendieck
~\cite{T1,T2,Mz,S2} asserts that the action of $G(K)$ on
$H^{2m}(X_a,\Ql)(m)$ is semisimple, i.e., the Frobenius
automorphism acts on $H^{2m}(X_a,\Ql)(m)$ as a semisimple linear
operator.

 The Tate conjecture is known to be true in certain cases, e.g., Fermat varieties satisfying certain numerical
 conditions ~\cite{T1,Sh,Y}; abelian varieties for $m=1$ \cite{T2}, certain classes of abelian
varieties with arbitrary $m$ ~\cite{Tn,Z2,Z3,LZ}, \cite{Milne2},
``almost all'' K3 surfaces ~\cite{ASD,N,NO}.

Let $Y$ be a {\sl cubic fourfold}, i.e., a smooth projective
hypersurface of degree $3$ in $\PP^5$ defined over $k$. It is
well-known that all odd $\ell$-adic Betti numbers of $Y_a$ do
vanish; it is also known that the second  and sixth Betti numbers
of $Y_a$ are equal to $1$ while its fourth Betti number is $23$
\cite{R}. Let $\Le \in H^2(Y_a,\Ql)(1)$ be the class of a
hyperplane section of $Y$. One may easily check that
$$H^2(Y_a,\Ql)(1)=\Ql\cdot\Le= H^2(Y_a,\Ql)(1)^{G(k)};
 H^6(Y_a,\Ql)(1)=\Ql\cdot\Le^3= H^6(Y_a,\Ql)(1)^{G(k)}.$$
It was proven by Rapoport \cite{R} that the Frobenius automorphism
acts on $H^{4}(X_a,\Ql)(2)$ as a semisimple linear operator.

 Recently N. Levin \cite{L} proved the Tate conjecture when
$X=Y$ is an {\sl ordinary} cubic fourfold. In this case the only
non-trivial case is $m=2$ and the theorem of Levin  asserts that
$H^4(Y_a,\Ql)(2)^{G(k)}$ is generated as $\Ql$-vector subspace by
algebraic classes of codimension $2$ (defined over $k$).

The aim of the present paper is to prove the Tate conjecture for
self-products $X=Y^r$ of an ordinary K3 cubic fourfold $Y$. For
example, if $r=2$ and $X=Y^2$ then the most interesting case is
$m=8$ and we prove, assuming $k$ ``sufficiently large", that
$H^8(X_a,\Ql)(2)^{G(k)}$ is generated as a $\Ql$-vector subspace
by algebraic classes of the following four types: 1) products of
codimension $2$ algebraic cycles on factors $Y$ of $Y^2$; 2) the
classes of graphs of Frobenius and its iterates; 3) the classes of
$Y \times \{\point\}$ and $\{\point \} \times Y$; 4)$\Le \times
\Le^3$ and $\Le^3 \times \Le$.
 The proofs are based on results and ideas of our previous papers ~\cite{Z4,Z5}
 where the Tate conjecture was proven for self-products of ordinary K3 surfaces .

The article is organized as follows. In Section \ref{finite} we
discuss tensor invariants of certain  $\ell$-adic
representations of $G(k)$.
  Section \ref{adic} treats cohomological $\ell$-adic representations.
 The Tate conjecture  is discussed in Section \ref{tate}.
  In Section \ref{prf} we will prove the Tate conjecture for powers of an ordinary
  cubic fourfold.

We write $\Z_+$ for the set $\{0,1,2,\ldots\}$ of non-negative
integers.
 Recall that $G(k)$ is procyclic and has a
canonical generator, namely, the arithmetic Frobenius automorphism
$$ \sigma_k : k(a) \rightarrow k(a), x  \rightarrow x^q$$
where $q$ is the number of elements of $k$. Clearly, $q$ is an integral power of $p$.
Another canonical generator of $G(k)$ is the geometric Frobenius automorphism $\varphi_k = \sigma_k ^{-1}$.

\section{Finite fields and  $\ell$-adic representations}
\label{finite}

Let $k$ be a finite field of characteristic $p$ consisting of $q$
 elements. We keep all notations of the Introduction connected with $k$. If
 $k'\subset k(a)$ is a finite overfield of $k$, then $k'(a)=k(a)$, the Galois group $G(k')$ of $k'$ is
 an open subgroup of finite index $[k':k]$ in $G(k)$ and $\varphi_{k'}=\varphi_k^{[k':k]}$.

 Let $\ell$ be a rational prime different from $p$. We
refer to ~\cite{Mn,T1,Z4} for the definition of the
one-dimensional $\Ql$-vector space $\Ql(1)$ and the corresponding
cyclotomic character
$$\chi_{\ell} : G(k) \rightarrow \Z_{\ell}^\times \subset \Q_{\ell}^{\times} =
\Aut\,\Q_{\ell}(1).$$ This character defines the Galois action on
$\Q_{\ell}(1)$. Notice, that
$$\chi_{\ell}(\sigma_k)=q, \quad \chi_{\ell}(\varphi_k)= q^{-1},.$$
We write $\Q_{\ell}(-1)$ for the  one-dimensional dual vector
space $\Hom_{\Q_{\ell}}(\Q_{\ell}(1), \Q_{\ell})$ with a natural
structure of the dual Galois module defined by the character
$\chi_{\ell}^{-1}$. To each integer $i$ one may attach a certain
one-dimensional $\Q_{\ell}$-vector space $\Q_{\ell}(i)$ with
Galois action defined by $\chi_l^i$. Namely,
$\Q_{\ell}(0)=\Q_{\ell}, \Q_{\ell}(i)=\Q_{\ell}(1)^{\otimes i}$ if
$i$ is positive and $ \Q_{\ell}(i)=\Q_{\ell}(-1)^{\otimes (-i)}$
if $i$ is negative.

For all integers $i,j$ there are natural isomorphisms of Galois modules
$$ \Ql(-i)=\Hom_{\Ql}(\Ql(i),\Ql),\; \Ql(i)\otimes_{\Ql}\Ql(j)=\Ql(i+j).$$

 Let
$$\rho : G(k) \rightarrow \Aut(V)$$
be an $\ell$-adic representation of $G(k)$, i.e., $V$ is a
finite-dimensional $\Q_{\ell}$-vector space and $\rho$ is a
continuous homomorphism \cite{S1}.  Clearly,
$$V^{G(k)}=V^{\rho(\varphi_k)}.$$
To each integer $i$ one may attach the \mbox{\rm   twisted}
$\ell$-adic representation
$$\rho[i] : G(k) \rightarrow \Aut(V(i))$$
where
$$V(i)=V\otimes_{\Ql}\Ql(i),\;
 \rho[i](\sigma)(v\otimes a)=\chi_l^i(\sigma)((\sigma v)\otimes a).$$

For example,
$$\rho[i](\sigma_k)(v\otimes a)=q^i(\rho(\sigma_k)(v)\otimes a),\;
\rho[i](\varphi_k)(v\otimes a)=q^{-i}(\rho(\varphi_k)(v)\otimes
a).$$ We have $V(0)=V, \Ql(i)(j)=\Ql(i+j).$ If
$\rho^*:G(K)\rightarrow \Aut(V^*) $ is the dual $\ell$-adic
representation then there are natural isomorphisms of Galois
modules $(V(i))^*=V^*(-i)$.

\begin{rem}
\label{twist}
 Let $\rho':G(k)\rightarrow \Aut(W)$ be (may be, another)
 $\ell$-adic representation of $G(K)$. Then we obtain natural isomorphisms of
 Galois modules
\begin{gather*}
 V(i)\otimes_{\Ql}W(j)=(V\otimes_{\Ql} W)(i+j),
\rho[i]\otimes \rho'[j]=(\rho\otimes\rho')(i+j),\\
\begin{split}
\Hom_{\Ql}(V,W) &= V^*\otimes_{\Q_{\ell}}W=V^*(-i)\otimes_{\Ql}W(i)\\
&= V(i)^*\otimes_{\Ql}W(i)=\Hom_{\Q_l} (V(i),W(i)).
\end{split}
\end{gather*}
\end{rem}
\begin{ex}
\label{twist1} If $u=\rho(\sigma)$ for some $\sigma\in G(k)$ then
one may easily check that
$$u_{[i]}=\chi_l^{-i}(\sigma)\rho [i](\sigma) \in \End_{\Q_l}(V(i)).$$
In particular, if $u=\rho(\sigma_k)$ then $u_{[i]}=q^{-i}\rho [i](\sigma_k)$. If
$u=\rho(\varphi_k)$ then $u_{[i]}=q^i\rho [i](\varphi_k)$.
\end{ex}

\begin{thm}[Theorem 3.1.4 of \cite{Z5}]
\label{th314}
  Let us assume that the $\ell$-adic
representation $\rho [i]$ is semisimple and consider the
characteristic
 polynomial
$$P_{\rho [i]}(t)=\Det(\id-t\rho [i](\varphi_k),V(i)).$$
Let  $R$ be the set of reciprocal roots of $P_{\rho [i]}(t)$.
Assume that either $1$ is the unique element of $R$, i.e.,
$R=\{1\}$ is the one-element set or
there exists a non-empty subset $B\subset R$ such that:
\begin{enumerate}
\item $B$ consists of multiplicatively independent elements; in particular,
 it does not contain 1
and does not meet $B^{-1}=\{\alpha^{-1}\mid \alpha \in B\};$
\item  Either
 $R$ coincides with the disjoint union of $B$ and $B^{-1}$ or $R$ coincides
 with the  disjoint union of $B$, $B^{-1}$ and the one-element set $\{1\}$.
\end{enumerate}
Then for each even natural number $2n$ all elements of
$(V(i)^{\otimes 2n})^{G(k)} $ can be presented as a linear
combination of  tensor products of $n$ elements of $(V(i)^{\otimes
2})^{G(k)}.$ Each element of $(V(i)^{\otimes (2n+1)})^{G(k)}$ can
be presented as a linear combination of tensor products of an
element of $V(i)^{G(k)}$ and $n$ elements of $(V(i)^{\otimes
2})^{G(k)}.$
\end{thm}

\begin{defn}
 Recall ~\cite[3.2]{Z4} that $\rho$ is called
{\sl semistable} if it enjoys the following property. If
$u=\rho(\sigma) \in \Aut(V)$ for some $\sigma\in G(k)$ and an
eigenvalue $\alpha$ of $u$ is a root of unity, then $\alpha=1.$ In
fact, in order to make sure that $\rho$ is semistable, it suffices
to inspect the eigenvalues only for $u=\rho(\varphi_k)$
~\cite[3.2.1]{Z4}.
\end{defn}

\begin{rem}
\label{ss}
 If $\rho$ is semistable and $k'$ is a finite overfield
of $k$ then the restriction of $\rho$ to $G(k')$ is also a
semistable $\ell$-adic representation of $G(k')$ and the
invariants of $G(k)$ and $G(k')$ coincide, i.e.,
$V^{G(k)}=V^{G(k')}.$

Conversely, for each $\rho$ there exists a positive integer $r$
such that if an eigenvalue $\alpha$ of $\rho(\varphi_k)$ is a root
of unity then $\alpha^r=1$. Now if $k_r\subset k(a)$ is the degree
$r$ extension of $k$ then every eigenvalue $\beta$ of
$\rho(\varphi_{k_r})=\rho(\varphi_k)^r$ that is a root of unity is
equal to $1$. This means that the restriction of $\rho$ to
$G(k_r)$ is semistable.
\end{rem}

\section{$\ell$-adic cohomology}
\label{adic} Let $Y$ be a smooth projective variety over $k$ and
$Y_a=Y\times k(a)$. Let $\ell$ be a rational prime $\neq p$. The
Galois group $G(k)$ acts on $Y_a=Y\times_k k(a)$ through the
second factor. For each non-negative integer $i$ this action
induces by functoriality the Galois action on the $i$th
$\ell$-adic \'etale cohomology group $H^i(Y_a,\Ql)$. We write
$\rho_{Y,i}$ for the corresponding \mbox{\rm   cohomological}
$l$-adic representation \cite{S2}
$$\rho_{Y,i}: G(k)\rightarrow \Aut(H^i(Y_a,\Q_{\ell})).$$
Let $F_{Y/k}:Y \rightarrow Y$ be the Frobenius endomorphism of the
$k$-scheme $Y$ \cite{T1}, \cite{Mn}. It is defined as the identity
map on points, together with the map $f\rightarrow f^q$ in the
structure sheaf. Let $\Frob_Y=F_{Y/k}\times \id_{k(a)}$ be the
corresponding $k(a)$-endomorphism of $Y_a$.
 We write $(\Frob_Y)_i$ for the endomorphism of $H^i(Y_a,\Ql)$ induced by $\Frob_Y$.

 \begin{rem}
 \label{ext}
If
 $k'\subset k(a)$ is a finite overfield of $k$ and
 $Y'=Y\times_k k'$ is the corresponding smooth projective variety over $k'$, then
 $$Y'_a=Y_a,\quad H^i(Y'_a,\Q_{\ell})=H^i(Y_a,\Q_{\ell})$$
 and $\rho_{Y',i}$ coincides with the restriction of $\rho_{Y,i}$
 to $G(k')$. In particular,
 $$\rho_{Y',i}(\varphi_{k'})=\rho_{Y,i}(\varphi_k)^{[k':k]}.$$
\end{rem}

\begin{rem}
\label{twist2}
 It is well-known ~\cite{T1,Mn} that  $$(\Frob_Y)_i=\rho_{Y,i}(\varphi_k).$$
 \end{rem}

{Let us consider the characteristic polynomial
$$P_{Y,i}(t) = \Det(\id - t \rho_{Y,i}(\varphi_k), H^i(Y_a,\Ql)).$$
A famous theorem of Deligne ~\cite{D2,D3} (conjectured by Weil)
asserts that $P_{Y,i}(t)$ lies in $1+t\Z[t]$, does not depend on
the choice of $l$ and all its (complex) reciprocal roots have
absolute values equal to $q^{i/2}$.

Notice, that in the case of cubic fourfolds this result was proven
earlier than the general case by Rapoport \cite{R} (inspired by
ideas of \cite{D1}). His paper also contain the proof of
semisimplicity of the action of Frobenius in the case of cubic
fourfolds. For Abelian varieties the semisimplicity was proven by
Weil (see \cite{Md}).

 Let $i=2m$ be an even non-negative integer. Let us consider the twisted cohomological $l$-adic representation
$$\rho_{Y,2m}[m]: G(k)\rightarrow \Aut(H^{2m}(Y_a,\Ql)(m)).$$
One may easily check (\cite{Z4}, 4.2) that
$$P_{Y,[m]}(t) = \Det(\id - t \rho_{Y,2m}[m](\varphi_k), H^{2m}(Y_a,\Ql)(m))=
P_{Y,2m}(t/q^m).$$
 Now, the theorem of Deligne implies that
 $P_{Y,[m]}(t)$ lies in $1+t\Z[1/q][t]$, does not depend on the choice of $l$
 and all its (complex) reciprocal roots have absolute values equal to 1.
 In other words, $P_{Y,[m]}(t)$ is a $q$-admissible polynomial in the sense
 of \cite{Z4}. In particular, if $\alpha$ is a (reciprocal) root of $P_{Y,[m]}(t)$
then $\alpha^{-1}$ is also one.

\begin{rem}
\label{extss} Clearly, there exists a positive integer $r$ such
for all roots $\alpha$ of $P_{Y,[m]}(t)$ we have $\alpha^r=1$ if
$\alpha$ is a root of unity. It follows from Remarks \ref{ss} and
\ref{ext} that if $k'$ is the degree $r$ extension $k_r$ of $k$
and $Y'=Y\times_k k'$ then $\rho_{Y',2m}[m]$ is semistable for all
$\ell \ne p$.
\end{rem}

 Let $d=\dim(Y)$. Let us consider the K\"unneth decomposition
$$H^{2d}(Y_a \times Y_a,\Ql)(d)=\oplus_{i=0}^{i=2d}H^{2d-i}(Y_a,\Ql)(d)
\otimes H^i(Y_a,\Ql)$$
of the $d$-twisted middle $l$-adic cohomology group of $Y_a \times Y_a=
(Y \times Y)_a$. Clearly, each element $c \in H^{2d}(Y_a \times Y_a,\Ql)(d)$ can be presented uniquely as a sum
$$c=\oplus_{i=0}^{i=2d} c_i \quad \mbox{\rm with }
 c_i\in H^{2d-i}(Y_a,\Ql)(d) \otimes H^i(Y_a,\Ql).$$
Notice that the K\"unneth decomposition is $G(k)$-equivariant. In particular,
$c \in H^{2d}(Y_a \times Y_a,\Ql)(d)$ is a $G(k)$-invariant if and only if all
$c_i$ are $G(k)$-invariant.



\begin{rem}
 Let $c \in H^{2d}(Y_a \times Y_a,\Ql)(d)$ be an
 {\sl algebraic cohomology class}, i.e., a linear combination of cohomology classes
 of closed irreducible codimension $d$ subvarieties on $Y_a \times Y_a$ \cite{T1}.
 It follows from results of \cite{KM} that all $c_i$ are also {\sl algebraic cohomology classes}.
\end{rem}

 Let $\mathfrak u$ be a $k(a)$-endomorphism of $Y(a)$. We write
$\Graph_{\mathfrak u}$ for the graph of $\mathfrak u$; it is a
$d$-dimensional irreducible closed subvariety of $Y_a \times Y_a$.
We write $\cl(\Graph_{\mathfrak u})$ for its
 $\ell$-adic cohomology class: it is an element of $H^{2d}(Y_a \times Y_a,\Ql)(d)$.
By functoriality, $\mathfrak u$ induces an endomorphism of
$H^i(Y_a,\Ql)$
 which will be denoted by
$${\mathfrak u}_i \in \End(H^i(Y_a,\Ql) ) =\End(V_i).$$

\begin{ex} If ${\mathfrak u}=\id$ then its graph is the diagonal $\Delta$
and therefore $\id_i$
is the identity endomorphism of $H^i(Y_a,\Ql)$.
\end{ex}

\begin{ex} If ${\mathfrak u}=\Frob_Y$ then according to \ref{twist2}
$${\mathfrak u}_i=(\Frob_Y)_i = \rho_{Y,i}(\varphi_k).$$
For each positive $j$ the $j$th power $\Frob_Y^j$ of $\Frob_Y$ is defined.
 As usual, if $j=0$, we put $\Frob_Y^j=\id$. Let us put
$${\mathfrak fr}_i =(\cl(\Graph_{\Frob_Y})_i).$$
It is known [7] that ${\mathfrak fr}_i$ can be presented as a
linear combination of $\cl(\Graph_{\Frob_Y^j})$ with rational
coefficients ($j\in\Z_+$). Notice, that it is well-known
~\cite{T1,Mn} that all $\cl((\Graph_{\Frob_Y})_i)$ are
 $G(k)$-invariants.
\end{ex}

\begin{rem} If $d=\dim(Y)=2m$ is even then there is the canonical
 isomorphism
$$H^d(Y_a,\Ql)(d) \otimes H^d(Y_a,\Ql)=
H^{2m}(Y_a,\Ql)(m) \otimes H^{2m}(Y_a,\Ql))(m).$$
\end{rem}

\begin{rem} If $d=\dim(Y)=2m$ is even then one may easily deduce from
\ref{twist}, \ref{twist1} and \ref{twist2}
 that
$q^m\rho_{Y,2m}[m](\varphi_k) = ((\Frob_Y)_i)_{[m]}$, i.e.,
$$g:=\rho_{Y,2m}[m](\varphi_k) = q^{-m}((\Frob_Y)_i)_{[m]}.$$
\end{rem}

\begin{thm}[Theorem 4.4 of \cite{Z5}]
\label{th44}
  Assume that \ $d=\dim(Y)=2m$ is even,
 $g=\rho_{Y,2m}[m](\varphi_k)$ is a semisimple linear operator and all its
eigenvalues different from 1 are simple.
 Then  the vector subspace of Galois invariants
$${(H^d(Y_a,\Ql)(d)} \otimes {H^d(Y_a,\Ql))^{G(k)}}
 \subset
{H^{2m}(Y_a,\Ql)(m)} \otimes {H^{2m}(Y_a,\Ql))(m)}$$
is generated by
${(H^{2m}(Y_a,\Ql)(m))^{G(k)}} \otimes {(H^{2m}(Y_a,\Ql))(m))^{G(k)}}$
 and all
$(cl(Graph_{\Frob_Y^j}))_d$
with $j\in\Z_+$.
 In particular,
${(H^d(Y_a,\Ql)(d)} \otimes {H^d(Y_a,\Ql))^{G(k)}}$
is contained in the vector subspace
 of
${H^{2d}(Y_a \times Y_a,\Ql)(d)^{G(k)}}$
generated by
 $${(H^{2m}(Y_a,\Ql)(m))^{G(k)}} \otimes {(H^{2m}(Y_a,\Ql))(m))^{G(k)}}$$
 and all
${\cl(\Graph_{\Frob_Y^j})}$ with $j\in\Z_+$.
\end{thm}

\section{The Tate conjecture}
\label{tate}
  We write $\Al^m(Y)$ for the $\Ql$-vector
subspace of $H^{2m}(Y_a,\Ql)(m))$ spanned by the cohomology
classes of all algebraic cycles of codimension $m$ on $Y$. It is
well-known [20] that
$$\Al^m(Y) \subset H^{2m}(Y_a,\Ql)(m))^{G(k)}.$$
Elements of
$H^{2m}(Y_a,\Ql)(m))^{G(k)}$ are called {\sl Tate classes} on  $Y$.

 Tate [19] conjectured that the following assertion holds true.
$$\Al^m(Y) = H^{2m}(Y_a,\Ql)(m))^{G(k)} .\leqno T(Y,m,k,l):$$

\begin{rem} Let $k'$ be a finite algebraic extension of $k$ and $Y':=Y\times_k k'$. It is known
[19] that if the assertion $T(Y',m,k',l)$ holds true then the
assertion  $T(Y,m,k,l)$ also holds true.
\end{rem}

 Recall \cite{Z4}, that $Y$ is called to be of K3 type
in dimension $2m$ if the characteristic polynomial $P_{Y,[m]}(t)$
is of K3 type, i.e. its $p$-adic Newton polygon \cite{K} enjoys
the following properties. There exists a non-zero rational number
$c$ such that the set of slopes is either $\{c,-c\}$ or
$\{c,-c,0\}$. In both cases slopes $c$ and $-c$ must have length
1.

For example, a K3 surface is of K3 type in dimension 2 if and only
if it is ordinary \cite{Z4}. An ordinary Abelian surface is of K3
type in dimension 2.

An ordinary cubic fourfold $Y$ is of K3 type in dimension 4.
Indeed, the Hodge numbers of a cubic fourfold (in dimension $4$)
are as follows \cite{R}.
$$h^{4,0}=h^{0,4}=0, h^{3,1}=h^{1,3}=1, h^{2,2}=21.$$
Since the Hodge polygon of an ordinary cubic fourfold coincides
with the Newton polygon, the $p$-adic Newton polygon of
$P_{Y,4}(t)$ admits the following description. Its set of slopes
is $\{1,2,3\}$; the length of both slopes $1$ and $3$ is $1$ while
the length of slope $2$ is $21$. This implies easily that the
Newton polygon of $P_{Y,[4]}(t)$ is of K3 type with the set of
slopes $\{-1,0,1\}$. (The length of its slopes $1$ and $-1$ is $1$
while the length of slope $0$ is $21$. )

\begin{rem} One may easily define motives of K3 type. The paper \cite{GY}
contains  examples of motives of K3 type arising from Fermat
varieties.
\end{rem}

\begin{thm}
\label{indep}  Let $d=\dim(Y)=2m$ is even. Assume that
 $\rho_{Y,2m}[m]$ is semistable.
Assume also $Y$ is of K3 type in dimension $2m$. We write $a(m,Y)$
for the multiplicity of $1$ viewed as a root of $P_{Y,[m]}(t)$.
Then
$$P_{Y,[m]}(t)=(1-t)^{a(m,Y)} P_{Y,\tr}(t)$$
where the polynomial $P_{Y,\tr}(t) \in \Q[t]$ enjoys the following
properties:
\begin{enumerate}
\item[(i)]
$P_{Y,\tr}(t)$ is irreducible over $\Q$;
\item[(ii)]
The the set $R_{Y,\tr}$ of  reciprocal roots of $P_{Y,\tr}(t)$
enjoys the following properties: if $\alpha \in R_{Y,\tr}$ then
$\alpha^{-1}$(= complex conjugate of $\alpha$) also belongs to
$R_{Y,\tr}$. In addition, $R_{Y,tr}$ does not contain roots of
unity.
\item[(iii)]
One may choose a subset $B \subset R_{Y,\tr}$ such that $B$ does
not meet $B^{-1}:=\{\alpha^{-1}\mid \alpha \in B\}$ and
$R_{Y,\tr}$ is the (disjoint) union of $B$ and $B^{-1}$.
\item[(iv)]
In addition, the set $B$ consists of multiplicatively independent
elements.
\end{enumerate}
\end{thm}

\begin{proof}
The semistability means that none of the roots of $P_{Y,\tr}(t)$
is a root of unity. Taking into account, that all the coefficients
of $P_{Y,\tr}(t)$ are real, we obtain (ii). Now the assertion
(iii) follows readily.

Recall that $P_{Y,[m]}(t)$ is of K3 type. It follows easily that
$P_{Y,\tr}(t)$ is also of K3 type. Recall that none of roots of
$P_{Y,\tr}(t)$ is a root of unity. Now the assertions (i) and (iv)
follow easily from general results about polynomials of K3 type
~\cite[th. 2.4.3 and 2.4.4]{Z4}.
\end{proof}

\begin{thm}[Theorem 5.3.1 of \cite{Z5}]
\label{531}
  Let $d=\dim(Y)=2m$ is even. Assume that
 $\rho_{Y,2m}[m]$ is semistable and semisimple.
If $Y$ is of K3 type in dimension $2m$  then ${(H^d(Y_a,\Ql)(d)}
\otimes {H^d(Y_a,\Ql))^{G(k)}}$ is generated as a vector subspace
of ${H^{2m}(Y_a,\Ql)(m)} \otimes {H^{2m}(Y_a,\Ql))(m)}$
 by
${(H^{2m}(Y_a,\Ql)(m))^{G(k)}} \otimes {(H^{2m}(Y_a,\Ql))(m))^{G(k)}}$ and all
$(\cl(\Graph_{\Frob_Y^j}))_d $ with $j\in \Z_+$. In particular,
${(H^d(Y_a,\Ql)(d)} \otimes {H^d(Y_a,\Ql))^{G(k)}}$ is contained in the vector
 subspace of ${H^{2d}(Y_a \times Y_a,\Ql)(n)^{G(k)}}$ generated by
$${(H^{2m}(Y_a,\Ql)(m))^{G(k)}} \otimes {(H^{2m}(Y_a,\Ql))(m))^{G(k)}}$$
 and all
$\cl(Graph_{Frob_Y^j}) $ with $j\in \Z_+$.
\end{thm}

\begin{cor}
\label{square}
 Let $Y$ be an ordinary cubic fourfold over a finite
field $k$ of characteristic $p$. Then there exists a finite
overfield $k'\supset k$  such that the ordinary cubic fourfold
$Y'= Y \times_k k'$ over $k'$  enjoys the following properties for
all primes $\ell \ne p$:

\begin{enumerate}
\item[(i)]
$\rho_{Y',4}[2]$ is semistable;
\item[(ii)]
${(H^4(Y'_a,\Ql)(2)} \otimes_{\Ql} {H^4(Y'_a,\Ql)(2))^{G(k')}}$ is
contained in the $\Ql$-vector
 subspace of ${H^{8}(Y'_a \times Y'_a,\Ql)(4)^{G(k')}}$ generated by
$\Al^2(Y') \otimes_{\Ql} \Al^2(Y')$ and all
$\cl(Graph_{Frob_Y{'}^j}) $ with $j\in \Z_+$.
\item[(iii)]
Let ${\mathfrak l} \in H^2(Y'_a,\Ql)(1)$ be the class of a
hyperplane section. Let ${\mathfrak a}$ be an effective 0-cycle on
$Y'$.
 (For instance, if $Y'(k')=Y(k')$ is non-empty then one may take as ${\mathfrak a}$ any $k'$-rational point on
 $Y$.)
Then $H^8(Y'_a \times Y'_a,\Ql)(4)^{G(k)}$ is generated as a
$\Ql$-vector subspace by  $\Al^2(Y'_a) \otimes_{\Ql} \Al_2(Y'_a)$,
${\mathfrak l} \otimes {\mathfrak l}^3$, ${\mathfrak l}^3 \otimes
{\mathfrak l}$, $\cl(\Graph_{\Frob_Y^j}) \; (j\in \Z_+)$
 and the classes of $Y' \times {\mathfrak a}$ and ${\mathfrak a} \times Y'$.
 In particular, the Tate conjecture holds true for for $Y'\times Y'$ and therefore also for $Y\times Y$.
 \end{enumerate}
\end{cor}

\begin{proof}
Let us choose a finite overfield $k'$ of $k$ such that all
$\rho_{Y',4}[2]$ are semistable: its existence follows from Remark
\ref{extss}.

It was already mentioned that, thanks to the theorem of Rapoport
\cite{R}, $\rho_{Y',4}[2]$ is semisimple. Since, thanks to the
theorem of Levine \cite{L}, $H^4(Y'_a,\Ql)(2))^{G(k')}=\Al^2(Y')$,
the  assertion (ii) follows from Theorem \ref{531}.

The assertion (iii) follows from (ii) combined with the
Galois-invariance of the K\"unneth decomposition for $H^8(Y'_a
\times Y'_a,\Ql)(2)$ and obvious equalities
$$H^8(Y'a,\Ql)(4)=H^8(Y'a,\Ql)(4)^{G(k')}=\Ql\cdot\cl({\mathfrak
a}),$$
$$H^2(Y'a,\Ql)(1)=H^2(Y'a,\Ql)(1)^{G(k')}=\Ql\cdot{\mathfrak
l},$$
$$H^6(Y'a,\Ql)(3)=H^6(Y'a,\Ql)(3)^{G(k')}=\Ql\cdot{\mathfrak
l}^3$$ (where $\cl({\mathfrak a})$ is the cohomology class of
${\mathfrak a}$).
\end{proof}

\begin{rem}
If $Y$ is an ordinary cubic fourfold over a finite field $k$ then
the already mentioned results of Levin  and Rapoport imply that
$$Al^2(Y)=\Num_2(Y)\otimes\Ql$$
where $\Num_2(Y)$ is the group of numerical equivalence classes of
cycles of codimension $2$ on $Y$. In particular, the rank of
$\Num_2(Y)$ coincides with the multiplicity of $1$ viewed as a
root of $P_{Y,[2]}(t)$ ~\cite[Th. 2.9, pp. 74--75]{T2}. Clearly,
if $\rho_{Y_1,4}[2]$ if semistable then the rank of $\Num_2(Y)$ is
$a(2,Y)=23- \deg(P_{Y,tr}(t))$.
\end{rem}

\begin{cor}
\label{product}
 Let $Y$ and $Z$ be  ordinary cubic fourfolds over a finite
field $k$ of characteristic $p$, enjoying the following
properties:

\begin{enumerate}
\item[(i)]
$\rho_{Y,4}[2]$ and $\rho_{Z,4}[2]$ are semistable;
\item[(ii)]
$a(2,Y)\ne a(2,Z)$.
\end{enumerate}
Let $\Le_Y\in H^2(Y_a,\Ql)(1)$ and $\Le_Z\in H^2(Z_a,\Ql)(1)$ are
classes of hyperplane sections of $Y$ and $Z$ respectively. Let
${\mathfrak a}_Y$ and ${\mathfrak a}_Z$ be effective zero cycles
on $Y$ and $Z$ respectively.

 Let us put $W=Y\times Z$. Then $(H^8(W_a,\Ql)(4))^{G(k)}$
is generated as a $\Ql$-vector subspace by  $\Al^2(Y_a)
\otimes_{\Ql} \Al_2(Z_a)$, $\Le_Y \otimes \Le_Z^3$, $\Le_Y^3
\otimes \Le_Z$
 and the classes of $Y \times {\mathfrak a}_Z$ and ${\mathfrak a}_Y \times Z$.
 In particular, the Tate conjecture holds true for for $Y\times
 Z$.
\end{cor}

 \begin{proof}
 It suffices to check that
 $$(H^4(Y_a,\Ql)(2)\otimes
 H^4(Y_a,\Ql)(2)))^{G(k)}=(H^4(Y_a,\Ql)(2))^{G(k)}\otimes
 (H^4(Y_a,\Ql)(2)))^{G(k)}.$$
 In order to do that it suffices to check that if $\alpha$ is a
 root of $P_{Y,[2]}(t)$ and $\beta$ is a
 root of $P_{Z,[2]}(t)$ then $\alpha\beta=1$ if and only if
 $$\alpha=1,\quad \beta=1.$$
 Let us prove it. Suppose $\alpha \ne 1$ and $\beta \ne 1$.
  Then $\alpha$ and $\beta^{-1}$ are roots
 of  $P_{Y,tr}(t)$ and $P_{Z,tr}$ respectively. But $P_{Y,tr}(t)$ and
 $P_{Z,tr}$ are $\Q$-irreducible polynomials with different
 degrees and therefore cannot have common roots. Hence
 $\alpha\ne\beta^{-1}$, i.e., $\alpha\beta\ne 1$.
 \end{proof}

 \begin{rem}
 Similar arguments  prove the Tate conjecture for
 $Y\times S$ where $S$ is an ordinary K3 surface over $k$ with
 semistable $\rho_{Y,2}[1]$ and $a(2,Y)\ne a(1,S)+1$. (The second
 Betti number of a K3 surface is $22$ while the fourth Betti number
 of a cubic fourfold is $23$.) Also by the same token one may
 prove the Tate conjecture for the product of two ordinary K3
 surfaces with different Picard numbers.
 \end{rem}

\section {Powers of fourfolds}
\label{prf}

\begin{thm}
\label{four}
  Let  $Y$ be a smooth geometrically irreducible $4$-dimensional
projective variety over $k$ such that the first and third Betti
numbers of $Y_a$ are zero and the second Betti number of $Y_a$ is
$1$. Let us assume that $\rho_{Y,4}[2]$
 is semisimple and semistable. Assume, in addition,  that
$Y$ is of K3 type in dimension $4$.
 Let $r>1$ be an integer and let us put
$X:=Y^r$. Then each cohomology class in
$(H^{2m}(X_a,\Ql)(m))^{G(k)}$ can be presented as a linear
combination of products of pullbacks of Tate
 classes on $Y$ and $Y^2$ with respect to the projection maps
 $X=Y^r \rightarrow Y, X=Y^r \rightarrow Y^2$.
 In particular, if for some prime $\ell$ the Tate conjecture holds true for $Y$ and $Y^2$
 then it is also true for $X$ with the same $\ell$.
\end{thm}

\begin{cor}
\label{cor1}
  Let $Y$ is an ordinary cubic fourfold over a finite field $k$.
  Let $r>1$ be an integer and let us put
$X:=Y^r$.
Then the Tate conjecture holds true for $X$.
\end{cor}

\begin{proof}[Proof of Corollary \ref{cor1}]
According to Corollary \ref{square} there exists a finite
overfield $k'\supset k$ such that if $Y'=Y\times k'$ then
$\rho_{Y',4}[2]$ is semistable and the Tate conjecture is true for
$Y'^2$. Recall \cite{L} that the Tate conjecture is valid for
$Y'$. Applying Theorem \ref{four} to $Y'/k'$, we conclude that the
Tate conjecture is valid for $Y'^r$. Since $Y'^r=(Y\times
k')^r=Y^r\times k'$, we conclude that the Tate conjecture is valid
for $Y^r=X$.
\end{proof}

\begin{proof}[Proof of Theorem \ref{four}]
Let $\Le \in H^2(Y_a,\Ql)$ be the class of a hyperplane section of
$Y$. Clearly,
$$H^8(Y'_a,\Ql)(4)=(H^8(Y'_a,\Ql)(4))^{G(k')}=\Ql\cdot\cl({\mathfrak
a}),$$
$$H^2(Y'_a,\Ql)(1)=(H^2(Y'_a,\Ql)(1))^{G(k')}=\Ql\cdot{\mathfrak
l},$$
$$H^6(Y'_a,\Ql)(3)=(H^6(Y'_a,\Ql)(3))^{G(k')}=\Ql\cdot{\mathfrak
l}^3$$ (where $\cl({\mathfrak a}$ is the cohomology class of an
effective $0$-cycle ${\mathfrak a}$) on $Y$. In particular, the
cohomology spaces $H^2(Y'_a,\Ql)(1),
H^6(Y'_a,\Ql)(3),H^8(Y'_a,\Ql)(4)$ consist of Tate classes.

We say that $c \in H^{2m}(X_a,\Ql)(m)$ is a {\bf decomposable}
 cohomology class if it can be presented as a linear combination of products
 of pullbacks of Tate classes on $Y$ and $Y^2$ with respect to the projection
 maps
$$X=Y^r \rightarrow Y, \; X=Y^r \rightarrow Y^2.$$
Clearly, linear combinations and $\cup-$products of decomposable cohomology classes are also decomposable one.

Let $r^{\prime} <r$ be a positive integer, $Y^r\rightarrow
Y^{r^{\prime}}$ any projection map. If $c \in
H^{2m}(Y^{r^{\prime}}_a,\Ql)(m)$ is a decomposable cohomology
class on $Y^{r^{\prime}}_a$ then its pullback is a decomposable
cohomology class in $H^{2m}(Y^r_a,\Ql)(m)=H^{2m}(X_a,\Ql)(m)$. If
$r=1$ or $r=2$ then each Tate class on $X_a=Y^r_a$ is decomposable
by obvious reasons.
    Clearly, in order to prove  Theorem \ref{four}, we have to check that  each Galois-invariant
    cohomology class in $H^{2m}(X_a,\Ql)(m)$ is decomposable.

 Let us look more thoroughly at  the cohomology of
$X_a=Y_a^r$. First, notice that the K\"unneth formula combined
with Poincar\'e duality \cite{Mn} implies that under our
assumptions all odd-dimensional cohomology groups of $X_a$ vanish.
In order to describe explicitly the even-dimensional cohomology
groups of $X_a$ let us fix a non-negative integer $m$ and consider
the set ${\mathcal M}(r,m)$ of maps
$$j: \{1,2 \ldots ,r\} \rightarrow \{0,1,2,3,4\} \quad \mbox{\rm with }
\sum_{i=1}^r j(i)=m.$$
Then the K\"unneth formula for $X_a=Y_a^r$ implies easily that
$$H^{2m}(X_a,\Ql)=
\sum_{j \in {\mathcal M}(r,m)}\otimes_{i=1}^r
H^{2j(i)}(Y_a,\Ql).$$

After the proper twist we obtain a canonical isomorphism of Galois modules
$$H^{2m}(X_a,\Ql)(m)=
\sum_{j \in {\mathcal M}(r,m)}\otimes_{i=1}^r
H^{2j(i)}(Y_a,\Ql)(j(i))$$ compatible with $\cup-$products. In
particular,
$$(H^{2m}(X_a,\Ql)(m))^{G(k)}=\sum_{j \in {\mathcal M}(r,m)} H_j^{G(k)} \quad
\mbox{\rm where }
H_j:=\otimes_{i=1}^r H^{2j(i)}(Y_a,\Ql)(j(i)).$$

 The symmetric group ${\mathbf S}_r$ of permutation in
$r$ letters acts on $X=Y^r$ in a natural way. By functoriality, it
acts on $H^{2m}(X_a,\Ql)(m)$ and this action commutes with Galois
action. Clearly, if
 $s \in {\mathbf S}_r$ and $c \in H^{2m}(X_a,\Ql)(m)$ then the cohomology class $c$ is decomposable if and only
 if $s^*c$ is decomposable.
Notice also that
$$s^*H_j=H_{js^{-1}}\; \forall j\in {\mathcal M}(r,m)$$
with the map $js^{-1}:\{1,2 \ldots ,r\} \rightarrow \{0,1,2,3,4\},
\quad
 js^{-1}(i):=j(s^{-1}(i))$.
Of course, the latter formula defines an obvious action of
${\mathbf S}_r $ on ${\mathcal M}(r,m)$. It follows that
$s^*(H_j^{G(k)})=(H_{js^{-1}})^{G(k)}$  for all $j\in {\mathcal
M}(r,m)$; in particular, all Galois-invariant cohomology classes
in $H_j$ are decomposable if and only if  all Galois-invariant
cohomology classes in $H_{js^{-1}}$ are decomposable. This implies
that in order to prove  Theorem \ref{four}, it suffices to check
that all Galois-invariant cohomology classes in $H_j$ are
decomposable for each {\bf non-decreasing} maps $j:\{1,2 \ldots
,r\} \rightarrow \{0,1,2,3,4\}$ from ${\mathcal M}(r,m)$.

So,  Theorem \ref{four} follows from the following assertion.

\begin{lem}
\label{633}
 Let $j$ be a non-decreasing map
$$j:\{1,2 \ldots ,r\} \rightarrow \{0,1,2,3,4\} \quad
\mbox{\rm with } \sum_{i=1}^r j(i)=m.$$ Let $H_j=\otimes_{i=1}^r
H^{2j(i)}(Y_a,\Ql)(j(i))$ be the correspondent K\"unneth chunk of
\linebreak $H^{2m}(X_a,\Ql)(m)=H^{2m}(Y_a^r,\Ql)(m).$ Then each
Galois-invariant cohomology class in $H_j$ is decomposable.
\end{lem}

\begin{proof}[Proof of Lemma \ref{633}]
We use induction by $r$. We already know that the Lemma is true
for  $r=1$ and $r=2$ . So, we may assume that $r>2$.

{\bf Case 1.} Assume that $j(1)<2$, i.e., $j(1)=0$ or $1$. Let us
consider the the projection map $\phi_1:X=Y^r\rightarrow Y$ on the
first  factor, the projection map $\phi:X=Y^r \rightarrow Y^{r-1}$
onto the product of last $(r-1)$ factors and a
 non-decreasing map
$$j^{\prime}:\{1,2 \ldots ,r-1\} \rightarrow \{0,1,2,3,4\}, \quad
j^{\prime}(i):=j(i+1)
\quad \mbox{\rm with } \sum_{i=1}^{r-1} j^{\prime}(i)=m-j(1).$$
Notice
 that $H^{2j(1)}(Y_a,\Ql)(j(1))$
consists of Tate classes and therefore
$\phi_1^*(H^{2j(1)}(Y_a,\Ql)(j(1))$ consists of decomposable
classes.
  Clearly,
\begin{gather*}
\begin{split} H_j &= H^{2j(1)}(Y_a,\Ql)(j(1))\otimes \;\otimes_{i=2}^r
H^{2j(i)}(Y_a,\Ql)(j(i)) \\
&= H^{2j(1)}(Y_a,\Ql)(j(1)) \otimes
H_{j^{\prime}}=\phi^*(H_{j^{\prime}}),
\end{split} \\
H_j^{G(k)}=H^{2j(1)}(Y_a,\Ql)(j(1))\otimes
(H_{j^{\prime}})^{G(k)}= \phi^*((H_{j^{\prime}})^{G(k)})
\end{gather*}
where $H_{j^{\prime}}=\otimes_{i=1}^{r-1}
H^{2j^{\prime}(i)}(Y_a,\Ql)(j^{\prime}(i))$ is a K\"unneth chunk
of $H^{2m-2j(1)}(Y_a^{r-1},\Ql)(m-j(1))$. By induction assumption,
all cohomology classes in $(H_{j^{\prime}})^{G(k)}$ are
decomposable and, therefore, their pullbacks with respect to
$\phi$ are also decomposable. Recall that
$\phi_1^*(H^{2j(1)}(Y_a,\Ql)(j(1))$ consists of decomposable
classes. This ends the proof of the Lemma in the case of $j(1)<2$.

{\bf Case 2.} Assume that $j(r)>2$, i.e. $j(r)=3$ or $4$. Let us
consider the projection map $\phi:X=Y^r \rightarrow Y^{r-1}$ onto
the product of first $(r-1)$ factors,   the projection map
$\phi_r:X=Y^r\rightarrow Y$ on the last factor and a
 non-decreasing map
$$j^{\prime}:\{1,2 \ldots ,r-1\} \rightarrow \{0,1,2,3,4\}, \quad
j^{\prime}(i):=j(i)
\quad \mbox{\rm with }
 \sum_{i=1}^{r-1} j^{\prime}(i)=m-j(r).$$
Notice that $H^{2j(r)}(Y_a,\Ql)(j(r))$
consists of Tate classes and therefore
$\phi_r^*(H^{2j(r)}(Y_a,\Ql)(2j(r)))$ consists of decomposable
classes. This implies that
\begin{gather*}
\begin{split} H_j &= \otimes_{i=1}^{r-1} H^{2j(i)}(Y_a,\Ql)(j(i))\:\otimes
H^{2j(i)}(Y_a,\Ql)(j(r)) \\
&= \phi^*(H_{j^{\prime}})\:\otimes
\phi_r^*(H^{2j(r)}(Y_a,\Ql)(2j(r))),
 \end{split} \\
H_j^{G(k)}= (H_{j^{\prime}})^{G(k)}\:\otimes
H^{2j(r)}(Y_a,\Ql)(4)= \phi^*((H_{j^{\prime}})^{G(k)})\:\otimes
\phi_r^*(H^{2j(r)}(Y_a,\Ql)(j(r)))
\end{gather*}
where $H_{j^{\prime}}=\otimes_{i=1}^{r-1}
H^{2j^{\prime}(i)}(Y_a,\Ql)(j^{\prime}(i))$ is a K\"unneth chunk
of $H^{2m-2j(r)}(Y_a^{r-1},\Ql)(m-j(r))$. By induction assumption,
all cohomology classes in $(H_{j^{\prime}})^{G(k)}$ are
decomposable and, therefore, their pullbacks with respect to
$\phi$ are also decomposable. Recall that
$\phi_r^*(H^{2j(r)}(Y_a,\Ql)(2j(r)))$ consists of decomposable
classes. This ends the proof of
 the Lemma in the case of $j(r)>2$.

{\bf Case 3.} Assume that $j(1)>1$ and $j(r)<3$. Since $j$ is
non-decreasing and may take on only values $0,1,2,3,4$, this
implies that $j(i)=2$ for all $i$ and, therefore,
$$H_j=\otimes_{i=1}^r H^{2j(i)}(Y_a,\Ql)(j(i))=
\otimes_{i=1}^r H^4(Y_a,\Ql)(2).$$ So, we have to prove that all
cohomology classes in
$$H_j^{G(k)}=
(\otimes_{i=1}^r H^4(Y_a,\Ql)(2))^{G(k)}$$ are decomposable.
Notice that semistability of $\rho_{Y,2}[1]$ implies that each
(reciprocal) root of $P_{Y,[1]}(t)$ which is a root of unity must
be equal to 1.  Now, the decomposabilty property of elements of
 $(\otimes_{i=1}^r H^2(Y_a,\Ql)(1))^{G(k)}$
 follows from the combination of Theorem \ref{indep} applied to $Y$ and $m=2$ and Theorem \ref{th314} applied to $V=H^4(Y_a,\Ql)$,
$\rho=\rho_{Y,2}$ and $i=2$.
 \end{proof}
  This ends the proof of Theorem \ref{four}.
  \end{proof}


\begin{thebibliography}{99}
\bibitem{Ar} M. Artin, {\em Supersingular K3 surfaces},
 Ann. Scient. \'Ecole Norm. Sup. 4e s\'erie {\bf 7} (1974), 543-568.

\bibitem{ArMz} M. Artin and B. Mazur, {\em Formal groups arising from algebraic
 varieties}, Ann. Scient. \'Ecole Norm. Sup. 4e s\'erie {\bf 10} (1977), 87-132.

\bibitem{ASD} M. Artin and H. P. F. Swinnerton-Dyer, {\em The Tate-Shafarevich
 conjecture for pencils of elliptic curves on K3 surfaces}, Invent. Math. {\bf 20} (1973), 279-296.

\bibitem{D1} P. Deligne, {\em La conjecture de Weil pour les surfaces K3}, Invent. Math. {\bf 15} (1972), 206-226.

\bibitem{D2} P. Deligne, {\em La conjecture de Weil}.I., Publ. Math. IHES
{\bf 43}(1974), 273-307.

\bibitem{D3} P. Deligne, {\em La conjecture de Weil}.II., Publ. Math. IHES
{\bf 52}(1980), 137-252.

\bibitem{GY} F. Q. Gouv\^ea  and N. Yui, {\em Observations on Fermat motives
 of K3 type}, J. Number Theory  {\bf 71} (1998), 203--226.


\bibitem{KM} N. Katz and W. Messing, {\em Some consequences of the Riemann
 hypothesis for varieties over finite fields}, Invent. Math. {\bf 23}(1974), 73-77.

\bibitem{K} N. Koblitz, $p$-{\em adic numbers}, $p$-{\em adic analysis and
 zeta functions}, Springer-Verlag, Berlin - New York, 1984.

\bibitem{LZ} H. W. Lenstra, Jr and Yu. G. Zarhin, {\em The Tate conjecture for
 almost ordinary Abelian varieties over finite fields},
Advances in Number Theory (F. Gouv\^ea and N. Yui, eds.),
Clarendon Press, Oxford, 1993, 179-194.

\bibitem{L} N. Levin, {\em The Tate conjecture for cubic
fourfolds over a finite fields}. Compositio Math. {\bf 127}
(2001), 1--21.

\bibitem{Mz} B. Mazur, {\em Eigenvalues of Frobenius acting on algebraic
 varieties over finite fields},  Proc. Sympos. Pure Math., vol. 29 ,
 Amer. Math. Soc., Providence, RI, 1975, 231--261.

\bibitem{Mn} J. S. Milne, {\em \'Etale cohomology}, Princeton University Press, 1980.

\bibitem{Milne2} J. S. Milne, The Tate conjecture for certain abelian varieties over finite fields.
Acta Arith. {\bf 100} (2001), no. 2, pp.135--166.


\bibitem{Md} D. Mumford, {\em Abelian varieties}, second edition, Oxford University Press, London, 1974.

\bibitem{N} N. O. Nygaard, {\em The Tate conjecture for ordinary K3 surfaces
 over finite fields}, Invent. Math. {\bf 74} (1983), 213--237.

\bibitem{NO} N. O. Nygaard and A. Ogus, {\em The Tate conjecture for  K3
 surfaces of finite height}, Ann. of Math. {\bf 122} (1985), 461--507.


\bibitem{R} M. Rapoport, {\em Complement \'a l'article de
P.Deligne} ``La conjecture de Weil pour les surfaces K3", Invent.
Math. {\bf 15} (1972), 227--236.

\bibitem{S1} J.-P. Serre, {\em Abelian} $\ell$-{\em adic representations and
 elliptic curves}, Second Edition, Addison-Wesley, 1989.

\bibitem{S2} J.-P. Serre, {\em Repr\'esentations} $\ell$-{\em adiques},
 Algebraic Number Theory (S. Yanaga, ed.),
 Japan Soc. Promotion  Sci., Tokyo, 1977, 177--193.

\bibitem{Sh} T. Shioda, {\em The Hodge conjecture and the Tate conjecture for
 Fermat varieties}, Proc. Japan Academy {\bf 55}(1979), 111--114.

\bibitem{Tn} S. G. Tankeev, {\em On cycles on Abelian varieties of prime
 dimension over finite or number fields}, Math. USSR Izvestija {\bf 22}(1984), 329--337.

\bibitem{T1} J. Tate, {\em Algebraic cycles and poles of zeta functions},
 Arithmetical Algebraic Geometry (O. F. G. Schilling, ed.), Harper and Row, New York, 1965, 93--110.

\bibitem{T2} J. Tate, {\em Endomorphisms of Abelian varieties over finite
 fields}, Invent. Math. {\bf 2}(1966), 134--144.

\bibitem{T3} J. Tate, {\em Conjectures on algebraic cycles in}
$\ell$-{\em adic cohomology}. In: Motives (U. Jannsen, S. Kleiman,
J.-P. Serre, eds.), Proc. Symp. Pure Math. {\bf 55}, Part 1
(1994), 71--83.


\bibitem{Y} N. Yui, {\em Special values of zeta functions of Fermat varieties
 over finite fields}, Number theory New York Seminar 1989-1990, Springer-Verlag, New York, 1991, 251--275.


\bibitem{Z2} Yu. G. Zarhin, {\em Abelian varieties of K3 type and}
 $\ell$-{\em adic representations},
 Algebraic Geometry and Analytic Geometry Tokyo 1990,
 ICM-90 Satellite Conference Proceedings, Springer-Verlag, Tokyo, 1991, 231--255.

\bibitem{Z3} Yu. G. Zarhin, {\em Abelian varieties of K3 type},
 S\'eminaire de Th\'eorie des Nombres, Paris 1990-91 (S. David, ed.),
 Progress in Math., vol. 108, , Birkh\"auser, 1993, 263-279.


\bibitem{Z4} Yu. G. Zarhin, {\em Transcendental cycles on ordinary K3 surfaces
 over finite fields}, Duke Math. J., {\bf 72}(1993), 65--83.

\bibitem{Z5} Yu. G. Zarhin, {\em  The Tate Conjecture for Powers of Ordinary K3 surfaces Over Finite
Fields}, J. Algebraic Geometry  {\bf 5} (1996), 151--172.


\end{thebibliography}
\end{document}